\documentclass[aop,preprint]{imsart}

\RequirePackage{amsthm,amsmath,amsfonts,amssymb}
\RequirePackage[numbers]{natbib}
\RequirePackage[colorlinks,citecolor=blue,urlcolor=blue]{hyperref}
\RequirePackage{graphicx}

\RequirePackage{mathrsfs}
\usepackage{mathrsfs}

\startlocaldefs

\theoremstyle{plain}

\newtheorem{theorem}{Theorem}[section]
\newtheorem{lemma}[theorem]{Lemma}

\theoremstyle{remark}

\allowdisplaybreaks[3]
\usepackage{color}
\newcommand{\bbA}{{\bf A}}
\newcommand{\bbP}{{\bf P}}
\newcommand{\bbE}{{\bf E}}

\newcommand{\bbC}{{\bf C}}
\newcommand{\bbD}{{\bf D}}
\newcommand{\bbS}{{\bf S}}
\newcommand{\bbR}{{\bf R}}
\newcommand{\bbT}{{\bf T}}

\newcommand{\bbI}{{\bf I}}

\newcommand{\bbB}{{\bf B}}
\newcommand{\bbX}{{\bf X}}

\newcommand{\bbY}{{\bf Y}}

\newcommand{\us}{{\underline s}}
\newcommand{\ug}{{\underline g}}
\newcommand{\uF}{{\underline F}}

\newcommand{\bqn}{\begin{eqnarray*}}
	\newcommand{\eqn}{\end{eqnarray*}}

\newcommand{\bqa}{\begin{eqnarray}}
	\newcommand{\eqa}{\end{eqnarray}}

\newtheorem{thm}{Theorem}[section]
\newtheorem{cor}[thm]{Corollary}

\endlocaldefs

\begin{document}
	
	\begin{frontmatter}
		\title{Exact Separation of Eigenvalues of Large Dimensional noncentral Sample Covariance Matrices}
		\runtitle{}
		
		\begin{aug}
			
			\author[A]{\fnms{Zhidong}~\snm{Bai}\ead[label=e1]{baizd@nenu.edu.cn}}
			\author[A]{\fnms{Jiang}~\snm{Hu}\ead[label=e2]{huj156@nenu.edu.cn}}
			\author[B]{\fnms{Jack W.}~\snm{Silverstein}\ead[label=e4]{jack@ncsu.edu}}
			,
			\and
			\author[A]{\fnms{Huanchao}~\snm{Zhou}\ead[label=e3]{zhouhc782@nenu.edu.cn}}

			\address[A]{KLASMOE and School of Mathematics and Statistics, Northeast Normal University, China \printead[presep={,\ }]{e1,e2,e3}}
			\address[B]{Department of Mathematics, North Carolina State University, USA \printead[presep={,\ }]{e4}}
			
		\end{aug}
		
		\begin{abstract}
			
			Let $ \bbB_n =(\bbR_n +\frac{1}{\sqrt n} \bbT^{1/2}_n \bbX_n)(\bbR_n +\frac{1}{\sqrt n} \bbT^{1/2}_n \bbX_n)^* $ where $ \bbX_n $ is a  $ p \times n $ matrix with independent standardized random variables, $ \bbR_n $ is a $ p \times n $ non-random matrix, representing the information,  and $ \bbT_{n} $ is a $ p \times p $ non-random nonnegative definite Hermitian matrix. 
			Under some conditions on $ \bbR_n \bbR_n^* $ and $ \bbT_n $, it has been proved that for any closed interval outside the support of the limit spectral distribution, with probability one there will be no eigenvalues falling in this interval for all $ p $ sufficiently large.
			The purpose of this paper is to carry on with the study of the support of the limit spectral distribution, and we show that there is an exact separation phenomenon: with probability one, the proper number of eigenvalues lie on either side of these intervals.
			
		\end{abstract}
		
		
		
	\end{frontmatter}

	\section{Introduction}
	
	Let $\bbB_n$ be a $ p \times p $ matrix defined by
	\begin{equation} \label{bb1} 
		\bbB_n =\frac{1}{n}(\bbR_n + \bbT^{1/2}_n \bbX_n)(\bbR_n + \bbT^{1/2}_n \bbX_n)^*,
	\end{equation}
	where $ \bbX_n=(x_{ij}) $ is a  $ p \times n $ matrix of independent and standardized random variables 
	$ (\bbE x_{ij} =0, \bbE \lvert x_{ij} \rvert^2=1  ) $, $ \bbR_n $ is a $ p \times n $ non-random matrix and $ \bbT_{n} $ is a $ p \times p $ non-random non-negative definite Hermitian matrix. 
	The model $\bbB_n$ can be viewed as the non-central sample covariance matrix and it is also referred in the literature (see \cite{zhou2022limiting}) to as the general information-plus-noise type matrix where the information is contained in the matrix $ (1/n)\bbR_n \bbR_n^* $ and the matrix $ \bbT^{1/2}_n \bbX_n $ is the additive noise.

	In this paper, we carry on with the study of the location of eigenvalues of the matrix $\bbB_n$ when $ p $ and $ n $ are large with the ratio $ y_n = \dfrac{p}{n} \to y \in (0,1) $. 
	Before that, there have been some results on the non-central matrix $\bbB_n$.
	Under certain conditions, especially assuming $\bbT_n$ and $\bbR_n \bbR_n^*$ are commutative, the result of the limiting spectral distribution (LSD) of $\bbB_n$ has been studied in terms of the empirical spectral distribution (ESD) function $F^{\bbB_n}$ in \cite{zhou2022limiting}. 
	And the analytic behavior of the LSD of the matrix $ \bbB_n $ have been studied in \cite{zhou2022Analysis}. It is shown that the LSD has a continuous derivative away from zero, the derivative being analytic wherever it is positive, and the determination criterion is available for its support.
	The support of the probability distribution plays a fundamental role in the study of the spectrum of $\bbB_n$. In \cite{bai2023no} , it proved that for any closed interval outside the support of the LSD of $\bbB_n$, with probability one there will be no eigenvalues falling in this interval for all $n$ sufficiently large.  
	The aim of this paper is to prove that the proper number of eigenvalues lie on either side of these intervals. The analog of he property was first proposed in \cite{bai1999exact} called "exact separation". 
	
	Our work has been originally motivated by the popular sample covariance matrix $\bbC_n= \frac{1}{n}\bbT^{1/2}_n \bbX_n \bbX_n^*\bbT^{1/2}_n  $ of $n$ samples of the random vector 
	$ \bbT^{1/2}_n X_{\cdot 1} $($ X_{\cdot 1}$ denoting the $j$-th column of $\bbX_n$). 
	Bai and Silverstein  established the complete results about the almost sure localization of the eigenvalues of $\bbC_n$ in the non-Gaussian case in \cite{bai1998no,bai1999exact}. 	It is shown that for a closed interval $[a,b]$ contained in an open interval $(c,d)$ out the supports of LSD's $F^{y_n,H_n}$ for almost all large $n$, the numbers of eigenvalues of $\bbC_n$ falling on either side of $[a,b]$ are exactly the same of eigenvalues of $\bbT_n$ falling on the corresponding sides of the interval $J \subset \mathbb{R}^+ $ which determines the interval $[a,b]$. 
	
	Similar ideas were also developed for the information-plus-noise matrix $$\bbD_n =\frac{1}{n}(\bbR_n + \sigma \bbX_n)(\bbR_n + \sigma \bbX_n)^*.$$ Using the same technique of \cite{silverstein1995strong, silverstein1995analysis}, Dozier and Silverstein established the LSD of $\bbD_n$ in \cite{dozier2007empirical} and derived the analytical behavior of the LSD of $\bbD_n$ in \cite{dozier2007analysis}.
	Such exact separation phenomenon was also exhibited for information-plus-noise matrix in \cite{loubaton2011almost}, they proved that a gap in the spectrum of $\bbD_n$ which corresponds a gap in the spectrum of $\bbR_n\bbR_n^*$ splits the spectrum of $\bbR_n\bbR_n^*$ exactly as that of $\bbD_n$ dealing with independent Gaussian variables $X_{ij}$, and later, Capitaine extended the results to the framework of non-Gaussian information-plus-noise type matrices investigated in\cite{capitaine2014exact}. In addition, they addressed the behaviour of the largest eigenvalues of the information-plus-noise spiked model. And both results were obtained when the information matrix $\bbR_n\bbR_n^*$ is low rank.
	
	The majority of known results are under the central condition, that is, the entries of $\bbX_n$ are zero mean (which can be extended to allowing the entries to have a common mean). 
	In fact, the large non-central random matrices also have significant implications that may be considered as an extension of non-central Wishart matrices, an important random matrix in multivariate linear regression under a non-null hypothesis.
	We turn now to the aim of this paper, that is, to extend the above results of exact separation for the  non-central random matrices $\bbB_n$ when $\bbR_n\bbR_n^*$ is full rank. Note that these results may also hold in the full rank case, not treated in this paper.

	The rest sections of this paper are organized as follows. In section 2, we review some important results of the general information-plus-noise type matrix $\bbB_n$ as well as some useful mathematical tools. Our main result and detailed proof are present in Section 3.

	\section{Preliminary Results}
	
	In this section, we recall some known results related to $\bbB_n$. For the convenience of readers, we will expand here to illustrate.
	To establish exact separation, it is necessary to review properties of the LSD of of $\bbB_n$.
	The best way in understanding both the LSD and analytic properties of the LSD is investigating the equations of the Stieltjes transforms of the LSD. Let $ F $ be any probability distribution, its Stieltjes transform is defined as 
	$$ s_{F}(z)=\int \dfrac{1}{\lambda-z} \mathrm{d}F(\lambda), 
	z\in \mathbb{C}^{+}\equiv \{ z\in \mathbb{C}: \Im z >0 \},  $$
	and $ F $ can be obtained by the inversion formula 
	\begin{align}
		F(b)-F(a)=\frac{1}{\pi}\lim_{v\rightarrow 0^{+}}\int_{a}^{b}\Im s_{F}(x+iv)\mathrm{d} x, \label{1.3}
	\end{align}
	where $ a $, $ b $ are continuity points of $ F $. 
	
	\subsection{Convergence of the ESD of $\bbB_n$}
	
	The main result of the LSD can be summarized in the following lemma (see \cite{zhou2022limiting} for more details).
	
	\begin{lemma}[Theorem 2.1 of \cite{zhou2022limiting}] \label{lemmaLSD}
		Under the assumptions:
		\begin{description}
			\item
			[(a)] $ \frac{p}{n}=y_{n}\rightarrow y $, as $\min\{p,n\}\to\infty$, and $ y\in(0,1] $.
			\item
			[(b)] The entries of  $ \bbX_{n}=(X_{ij}) $ are  independent and normalized, i.e., with zero mean and unit variance. Also, satisfy the Lindeberg condition: for any $\eta>0$,
			\begin{equation*}
				\frac1{\eta^2 pn}\sum_{i\le p; j\le n} E\lvert X_{ij}\rvert^2I_{\{\lvert X_{ij}\rvert >\eta \sqrt{n}\}}\to 0.
			\end{equation*}
			\item
			[(c)] $ \bbR_{n}\bbR_{n}^{*} $ and $ \bbT_{n} $ are commutative.
			\item
			[(d)] As $\min\{p,n\}\rightarrow\infty$, the two-dimensional distribution function 
			$ H_n(u,t) = p^{-1} \sum_{i=1}^{p} I(u_i \le u, t_i \le t) $ converges weakly to a nonrandom limit distribution $ H(u, t) $, where $ u_{i}, t_{i} $ are the paired eigenvalues of 
			$ \frac{1}{n}\bbR_{n}\bbR_{n}^{*} $ and $ \bbT_{n} $,  respectively.	
		\end{description}
		
		Then, with probability one, $ F^{\bbB_n} $ converges in distribution to $ F, $ a nonrandom probability distribution function, whose Stieltjes transform $ s = s_F(z) $ satisfies the equation system
		\begin{equation}
			\left\{
			\begin{aligned}
				s=\int\frac{\mathrm{d}H(u,t)}{\frac{u}{1+yg}-(1+yst)z+t(1-y)}, \\
				g=\int\frac{t\mathrm{d}H(u,t)}{\frac{u}{1+yg}-(1+yst)z+t(1-y)}. \label{712.1}
			\end{aligned}
			\right.
		\end{equation}	
		Moreover, for each $ z\in \mathbb{C}^{+} $, $ (s, g) $ is the unique solution to $ \eqref{712.1} $ in  $ \mathbb{C}^{+} $.	
	\end{lemma}
	
	For convenience, we also consider the convergence of the ESD of $$\underline{\bbB}_n=\frac{1}{n}(\bbR_n+\bbT^{\frac{1}{2}}_n \bbX_n)^*(\bbR_n +  \bbT^{\frac{1}{2}}_n \bbX_n). $$
	The eigenvalues of the matrix $ \underline{\bbB}_n $ are the same as those of the matrix $ \bbB_n $ except $ \lvert n-p \rvert $ zero eigenvalues. 
	Let $F^{\underline{\bbB}_n}$ denote the almost sure limiting distribution function of the eigenvalues of $\underline B_n$. 
	Therefore, it is easy to verify the  following relations
	\begin{align*}
		F^{\underline{\bbB}_n}=\left(1-\frac{p}{n} \right)I_{[0,\infty]}+\frac{p}{n} F^{ \bbB_n},
	\end{align*} 
	where $I_A$ denoting the indicator function of the set $A$. 
	Then making a variable transformation
	\begin{align*}
		&\underline{s}(z)=-\dfrac{1-y}{z}+ys(z),  \\ 
		&\underline{g}(z)=-\dfrac{1}{z(1+yg(z))},   
	\end{align*}
	where $ \underline{s}(z)$ is the Stieltjes transform of the LSD 
	$ F^{\underline{\bbB}_n} $ and $ \underline{s}_{n}(z)=s_{F^{\underline{\bbB}_n}}(z) $ and  
	$ \underline{g}(z)$ is the limit of  $ \underline{g}_n=\frac{1}{p}\mathrm{tr} \bbT_n \left(\underline{\bbB}_n-z\bbI \right)^{-1} $.
	Because Stieltjes transform is invertible, then the equations in \eqref{712.1} become
	\begin{equation} 
		\begin{aligned}
			z=-\dfrac{1-y}{\underline{s}}-\frac{y}{\underline{s}}\int \dfrac{ \mathrm{d} H(u,t)}{1 + u\underline{g}(z) + t\underline{s}(z)},
			\\
			z=-\dfrac{1}{\underline{g}}+y\int \dfrac{t \mathrm{d} H(u,t)}{1 + u\underline{g}(z) + t\underline{s}(z)}. \label{713}
		\end{aligned}
	\end{equation}	
	Equating the two equations in \eqref{713}, we get 
	$$ -\frac{1-y}{\underline{s}} - \frac{y}{\underline{s}}\int \dfrac{\mathrm{d}H(u, t)}{1 + u\underline{g} + t\underline{s}} = -\frac{1}{\underline{g}} + y \int \dfrac{t\mathrm{d}H(u,t) }{1+u\underline{g}+ t\underline{s}}.  $$
	This is equivalent to
	\begin{equation}\label{7.3.2} 
		y\ug^2\int \dfrac{u\mathrm{d}H(u,t)}{1+u\underline{g}+t\underline{s}}+{\underline{s}} -{\underline{g}}=0.  
	\end{equation}	In addition, we consider several important quantities for the subsequent proof.
	Write $\underline{s}=\underline{s}_1+i\underline{s}_2 $, $ \underline{g} = \underline{g}_1+i\underline{g}_2 $, $ z\underline{s} = (z\underline{s})_1 + i(z\underline{s})_2 $, and $z\underline{g} = (z\underline{g})_1 + i(z\underline{g})_2$. 
	Fix $z=x+iv\in\mathbb{C}^+$. 
	Multiplying by $\underline s$ on both sides to \eqref{713} and comparing the imaginary part of the resulting equation, we obtain
		$$  (z\underline{s})_2 = cA_1 \underline{g}_2 + cB_1 \underline{s}_2,  $$
		$$  v =\frac{\underline{g}_2}{\vert \underline{g}\vert^2}-(c A_2 \underline{g}_2 + cB_2\underline{s}_2),  $$
		where
		\begin{align}
			\begin{aligned}
				A_j&=&\int \dfrac{ut^{j-1}\textrm{d} H(u,t)}{\vert 1 + u\underline{g} + t\underline{s}\vert^2}, j=1,2,
				\\
				B_j&=&\int \dfrac{t^{j}\textrm{d} H(u,t)}{\vert 1 + u\underline{g} + t\underline{s}\vert^2}, j=0,1,2. \label{7216}
			\end{aligned}
		\end{align}	
	 Lemma 1 of \cite{zhou2022Analysis} has been proved that all the four quantities $ A_j $, $ B_j $, $ j=1,2 $ are bounded.

	\subsection{Properties of the LSD of $\bbB_n$ and characterization of its support }	
	
	The behaviour of the Stieltjes transform of the LSD of $\bbB_n$ around the real axis is fundamental to evaluate its support of $F$. It is necessary at this point to review some of the properties of $ F $ and $ s_F $.  From \eqref{712.1}, much of the analytic behavior of $ F $ can be inferred in \cite{zhou2022Analysis}. 	And we recall the main properties in the following lemma. 

	\begin{lemma}[Theorem 3 of \cite{zhou2022Analysis} ]\label{lemma3} 
		Assume $c\leq1$ and the conditions imposed in the limiting $H$ in Lemma \ref{lemmaLSD}.  
		Suppose $\us(x) $ is the solution to \eqref{713} for $x\neq0$. 
		Then $ x \in \mathbb{R}\backslash\{0\} $, $ \lim_{z\in \mathbb{C}^+ \to x}s_{\underline F} (z) \equiv \us(x) $ exists. 
		The function $\us$ is continuous on $ \mathbb{R}\backslash\{0\} $, and $ F $ has a continuous derivative $ f $ on $ \mathbb{R}\backslash\{0\} $ given by $ f(x) = \frac{1}{\pi} \Im \us(x) $. Furthermore, if $ \Im \us(x) > 0$  $ (f(x) > 0) $ for $ x \in \mathbb{R}^+ $, then the density $ f $ is analytic about $ x $.
	\end{lemma}

	 Moreover, most important for this paper is to understanding the support of $F$. 
	 For any probability distribution function $ G $, let $ S_G $ and $ S_G^c $ denote the support of $ G $ and the complement of its support, respectively. By definition of support, we have $ S_F \subset [0, \infty)  $ and  $ S_H \subset [0, \infty)\times (0, \infty) $.
	 It is shown in \cite{zhou2022Analysis} that on any closed intervals outside the support of $ F $, 
	 $ s_F(x) $ exists and is increasing. Therefore, the inverse $ \eqref{1.3} $ can be used to identify these intervals. 
	 The following lemma has been proved in \cite{zhou2022Analysis}.
	 	\begin{lemma}[Theorem 5 of \cite{zhou2022Analysis} ]\label{lemma2}
	 	Assume $c\leq1$ and the conditions imposed on $H$ in  Lemma \ref{lemmaLSD} hold.  
	 	Let $x_0 \in S^c_{\underline{F}}\cap\mathbb R^+$,
	 	\begin{description}
	 		\item[(a)] Then $ \underline{s}(z)=\int (t-z)^{-1}\mathrm{d}\underline{F}(t) $ is analytic in a neighborhood $ D_{x_0} $ of $ x_0 $ and there exists a co-solution $ \underline{g}(z) $ which is also analytic in $ D_{x_0} $. 
	 		The triple $ (x, \underline{s}(x), \underline{g}(x)) $, $x\in D_{x_0}\cap\mathbb{R}^+$ is an extended solution to  
	 		\eqref{713} with $V < x$.
	 		\item[(b)] For any support point $ (u, t) $ of $ H $, $ u\underline{g}(x_0) + t\underline{s}(x_0)  \neq -1 $.
	 	\end{description}
	 	On the other hand, if $x_0,\us_0,\ug_0$, with $x_0>0$,  form a real extended solution to \eqref{713} satisfying (b), then from \eqref{7.3.2}, there exists a real analytic function $x=x(\ug)$, defined in an interval containing $\ug_0$ which satisfy \eqref{713}, and if $x'(\ug_0) \neq0$, then $x_0\in S_{\uF}^c$.
	 \end{lemma}

	 As for whether $F$ places any mass at $0$, it is also shown in \cite{zhou2022Analysis} that when 
	 $y\leq 1$, the LSD $F$ determined by \eqref{712.1} has no mass at zero.

	\subsection{No eigenvalues outside the support of the LSD of $\bbB_n$}
	
	The focus of \cite{bai2023no} is on intervals   $[a,b]\subset(c,d)$ lying outside the union of supports of $ F^{c,H} $ and $F^{c_n,H_n}$, for all large $n$, it proves that there is no eigenvalues of $\bbB_n$ falling in $a,b]$, where $F^{c,H}$ is the LSD of $\bbB_n$ and $F^{c_n,H_n}$ is an $F^{c,H}$ with $c$ and $H$ replaced by $c_n=p/n$ and $H_n$.
	
Roughly speaking, 
	when $ n $ is large, with probability one, there are no sample eigenvalues of $\bbB_n$ falling into the limiting spectral gaps.  
	
	More exactly,  the following result  for $\bbB_n$ is established  in\cite{bai2023no}.
	
	\begin{lemma}[Theorem 1.1 of \cite{bai2023no}]\label{lemma1}
		Assume that
		\begin{enumerate}
			\item[(a)] $ [a, b] \subset (c, d) \subset S^c_{\underline{F}^{y_n,H_n} }  $, with $ c > 0 $ for all large $ n $;
			\item[(b)] The matrix $\bbX_n$ is the $p\times n$ upper-left conner of the double array of random variables $ x_{ij} $ having means zero, variances one, second moments zero if complex and  there  is a random variable 
			$ X $ with finite fourth moment such that for a constant $ K $ and for all 
			$ x > 0 $ 
			$$\frac{1}{p} \sum_{i =1}^{p} \bbP(\lvert x_{ij}\rvert> x) \leq K\bbP(\lvert X \rvert) > x) , $$ and $$ \frac{1}{n}  \sum_{j=1} ^{n}\bbP(\lvert x_{ij} \rvert > x ) \leq K\bbP(\lvert X \rvert) > x) ; $$
			
			\item[(c)] There exists a positive function $ \psi(x) \uparrow \infty $ as $ x \to \infty$, and $ M > 0  $ such that
			$$\max_{ij} \bbE \lvert x^2_{ij} \rvert \psi(\lvert x_{ij} \rvert) \leq M; $$
			\item[(d)] $ n = n(p) $ with $ y_n = {p}/{n} \to y > 0 $  as $ n \to \infty $;
			\item[(e)] For $ n = 1, 2, \dots, $  $\bbR_n $ is a $ p \times n $ nonrandom matrix with $\dfrac{1}{\sqrt{n} } \bbR_n  $ uniformly bounded in special norm for all $ n $;
			\item[(f)] The matrix $ \bbT_n $ is uniformly bounded in spectral norm and 
			$ \lambda_{-1} \leq K $ for some constant $ K $, and is also commutative with $ (1/n)\bbR_n\bbR^*_n$ and their joint spectral distribution $ H_n(u, t)  $ tends to a proper distribution $ H(u, t) $, where $ \lambda_{-1} =\int t^{-1} \mathrm{d} H(u,t) $. 
		\end{enumerate}
		
		Then, we have that
		\begin{align} \label{7442}
			\bbP(\text{no eigenvalues of } \bbB_n \ \text{appear in} \ [a, b] \ \text{for all large} \ n) = 1. 
		\end{align}	
	\end{lemma}

		As mentioned in \cite{bai2023no}, assumptions (b)-(c) allow for the $ x_{ij} $ to depart from merely being i.i.d.. After suitable truncation, centralization, and scaling of the $ x_{ij} $'s one can assume these variables to be uniformly bounded.


\subsection{Other useful lemmas}

\begin{lemma}\label{lemma5.11}
	Assume that the entries of $ \{x_{ij} \} $ are a double array of independent complex random variables with mean zero, variance $ \sigma^2 $, and satisfy the assumptions (b) -- (e) of Lemma \ref{lemma1}. 
	Let $ \bbX_n = (x_{ij}; i \leq p, j \leq n) $ be the $ p \times n $ matrix of the upper-left corner of the double array. 
	Then, with probability one, we have
		$$  -2\sqrt{y}\sigma^2 \leq  {\lim\inf}_{n\to \infty} \lambda_{\min}(\bbS_n - \sigma^2(1 + y)\bbI_n)
		\leq {\lim \inf}_{n\to \infty} \lambda_{\max}(\bbS_n -\sigma^2(1 + y)\bbI_n)\leq 2\sqrt{y}\sigma^2,   $$
	where $\bbS_n=n^{-1}\bbX_n\bbX_n^*$.
	
\end{lemma}

The proof of Lemma \ref{lemma5.11} is proved in \cite{bai2023no} under the truncation and centralization, so the detailed
proof of the lemma is omitted. For details, the readers are referred to \cite{bai2023no}.

\begin{lemma}[Theorem A.46. of \cite{bai2010spectral} ]\label{lemma5}
	Let $\bbA$ and $\bbB$ be two $n \times p$ complex matrices. Then,
	$$\max_k \vert s_k(\bbA) - s_k(\bbB)\vert \leq \| \bbA - \bbB \|,$$ where $s_k(\bbA)$ is the $k$th largest singular value of the matrix $\bbA$.
	If $\bbA$ and $\bbB$ are Hermitian, then the singular values can be replaced by eigenvalues; i.e.,
	$$\max_k \vert \lambda_k(\bbA) - \lambda_k(\bbB)\vert \leq \| \bbA - \bbB \|.$$
\end{lemma}

	\section{Main Result}	
	
	In this section, we will study the exact separation problem of $\bbB_n$. For better characterize the exact separation, we equivalently reformulate the model $\bbB_n$ as
	\begin{align} \label{751}
		\bbB_n =(\bbR_n + \frac{1}{\sqrt{n}}\bbT^{1/2}_n \bbX_n)(\bbR_n + \frac{1}{\sqrt{n}}\bbT^{1/2}_n \bbX_n)^*.
	\end{align}
	The change is only redefined $\frac1{\sqrt n} \bbR_n$ as $\bbR_n$. It is not difficult to see that the LSD of $\bbB_n$ exists with probability one and its Stieltjes transform is given by \eqref{712.1} or equivalently by \eqref{713}. 
	
	When $\bbT_n =\sigma^2 \bbI_n$, we  recall the results of \cite{couillet2011deterministic,capitaine2014exact} that with probability one, the number of
	eigenvalues of $\bbB_n$ and $\bbR_n \bbR_n^* $ lying on one side of their respective intervals
	are identical for all large $n$. 
	Then considering the model \eqref{751}, we suppose that $ [a, b] $ is the interval defined in Lemma \ref{lemma1}, and based on the conditions about $ \bbR_n $, $ \bbT_n $, define
	$$ h_j(x) = u_j \underline{g}(x) + t_j \underline{s}(x), x \in [a, b], $$
	 where $ u_j $, $ t_j $ are the paired eigenvalues of  $\bbR_n \bbR_n^*$ and $\bbT_n$, respectively.
	
	Then, we have the following theorem.
	\begin{theorem}\label{th731}
		Under the assumptions of Lemma \ref{lemma1}, with probability one,
		when $ n $ is large, for each $ j \leq n $, the function
		$ h_j(x) $ is either $ > -1 $ for all $ x \in [a, b] $ or $ < -1 $ for all $ x \in [a, b] $. 
		Also, the number of eigenvalues of $ \bbB_n $ that are below $ a $ is exactly equal the number of 
		$ h_j(x) $ is below $ -1 $; and the number of eigenvalues of $ \bbB_n $ that are above $ b $ is exactly equal the number of $ h_j(x) $ is above $ -1 $.
	\end{theorem}
	
	Before starting the proof of the theorem, we first truncate and renormalize the random variables of  $\bbX_n$. We define $ \hat{\bbX}= (\hat{X}_{ij} )_{p\times n} $ and $$ \hat{X}_{ij}=X_{ij} I (\lvert X_{ij}\rvert <C) -\bbE X_{ij} I(\lvert X_{ij}\rvert <C) $$ for some constant $ C $,
	and define 
	$$ \hat{\bbB}_n=(\bbR_n+\frac{1}{\sqrt{n}}\bbT^{1/2}_n\hat{\bbX})(\bbR_n +\frac{1}{\sqrt{n}}\bbT^{1/2}_n \hat{\bbX})^*.  $$
	By Lemma \ref{lemma5.11} and Corollary 7.3.8 of \cite{hom1985matrix}, we have
	\begin{align*}
		\max_{i\leq p} \lvert \lambda_i (\bbB_n) - \lambda_i(\hat{\bbB}_n)\rvert  \leq  \dfrac{1}{\sqrt{n}}\lVert \bbT^{1/2}_n (\bbX - \hat{\bbX}) \rVert  
		\leq \sqrt{\lVert \bbT_n \rVert} (1 + \sqrt{y}) \sqrt{\bbE \lvert X_{ij}^2 \rvert I( \lvert X_{ij} \rvert  \geq  C)},
	\end{align*}
	which can be arbitrarily small when $ C $ is large.
	
	Another truncations of $\bbR_n $ and $\bbT_n$ can be found in \cite{zhou2022Analysis}, which ensured that $\bbR_n \bbR_n^* $ and $\bbT_n$ are of bounded norm for all $n$, and the details are not covered here.

\subsection{Spectral Gap Dependence on $ y $ as Tending to $ 0 $}	
	
	Now, using the same approach as in \cite{bai1999exact}, we set up a series of new models such that 
	$ y = y_\ell = p/n_\ell \to 0 $ as $ \ell \to \infty $. 
	To this end, we have to consider the model as 
	$$ \bbB_{n,\ell} = (\bbR_{n_\ell} +n^{-1/2}_{\ell} \bbT_n^{1/2} \bbX_{n_\ell})(\bbR_{n_\ell} +n^{-1/2}_{\ell} \bbT_n^{1/2} \bbX_{n_\ell})^*,  $$
	where $ n_{\ell} = n+\ell M $, and $ M/n \to \tau > 0 $ a small positive number. 
	That means, for each model $\bbB_{n,\ell}$, we extend the $\bbX_{n_{\ell}}$ with $M$ columns to obtain $\bbX_{n_{\ell+1}}$ and construct $\bbB_{n,\ell+1}$ accordingly.
	Also, we add $0$ columns to $\bbR_{n_{\ell}}$ to match the order of $\bbB_{n_{\ell}}$, enabling feasible matrix addition.
	
	Suppose that $[a, b]$ represents the spectral gap for $\bbB_n = \bbB_{n,0}$. 
	We will denote a series $[a_{\ell}, b_{\ell}]$ as the spectral gaps for $\bbB_{n,\ell}$ at each $\ell$ and demonstrate that the number of eigenvalues of $\bbB_{n,\ell}$ greater than $b_{\ell}$ (or smaller than $a_{\ell}$) is the same as the number of eigenvalues of $\bbB_{n,\ell+1}$ greater than $b_{\ell+1}$ (or smaller than $a_{\ell+1}$, correspondingly).
	To achieve this, we need to define $[a_{\ell}, b_{\ell}]$ and prove that $b_{\ell}-a_{\ell}$ increases as $\ell\to\infty$.

	We will consider solutions for \eqref{713} in the following way: for each $\underline{g}\in \mathbb{C}^+ $
	and $ y\leq 1 $, define $\us$ by \eqref{7.3.2} which exists uniquely and belongs to $ \mathbb{C}^+ $.
	Further, by the second equation in \eqref{713} we may find $ z = z_{y,H(\ug,\us,y)}\in \mathbb{C}^+ $ and get a set of solutions to \eqref{713}. 
	According to Lemma \ref{lemma3}, we can extend the solutions to \eqref{713} to the real axis. 
	We remind the readers that for any real $\ug$ and $y$, there is a unique real $\us$ solves \eqref{7.3.2} and consequently, there is a unique real $z$ which solves the second equation of \eqref{713}. 
	The triple $ (\ug, \us, x) $ is a set of solutions to \eqref{713} but may not be a set of extended solutions to it, that means, the limit of $ \us(z) $ when $ z\stackrel{\mathbb{C}^+}{\to}x + 0i$ may not be $ \us(x) $. According to Lemma \ref{lemma2}, this limit holds if and only if $ x \in S^c(F^{y,H}). $
	Therefore, we may use the positiveness of the derivative of $ x $ with respect to
	$ \ug $ to examine $x$ is outside the support of $F^{y,H}$.
	
	When the boundary point $ x $ of the support of $F^{y,H}$ is infinity, then by \eqref{713}, the corresponding $ \us $ and $ \ug $ are both $ 0 $. 
	Hence, for smaller $ y $, $ (0, 0, \infty) $ is still a set of solution. 
	Hence, $ \infty $ is still a boundary point of the support of $F^{y,H}$. 
	Consequently, when $ y $ decreases, the upper boundary of the spectral gap doesn't become smaller. Therefore, we only need to consider the case where the upper boundary of the spectral gap is finite.

	Recall the definition of \eqref{7216}, 
	\begin{align*}
		\begin{aligned}
			A_j&=&\int \dfrac{ut^{j-1}\textrm{d} H(u,t)}{\vert 1 + u\underline{g} + t\underline{s}\vert^2}, j=1,2,
			\\
			B_j&=&\int \dfrac{t^{j}\textrm{d} H(u,t)}{\vert 1 + u\underline{g} + t\underline{s}\vert^2}, j=0,1,2. 
		\end{aligned}
	\end{align*}	
	by conclusion (b) of Lemma 1 of \cite{zhou2022Analysis}, it has been proved that all the five quantities $ A_i $, $ B_j $ are all bounded for all $ z\in\mathbb{C}^+ $. Hence, all $ x $ are in the complement of the support of $F^{y,H}$ by Fatou Lemma. 
	That means $ (1+u\ug + t\us)^{-2}$ and its product with $t$, $ t^2 $, $ ut $ are integrable with respect to $H $, and thus the derivative of $ x $ with respect to $ \ug $ exists and is continuous. 
	Hence, the boundary points of spectral gaps are zero.
	
	A spectral gap is an open interval $ (c, d) $ that is an interval of the complement of support of $F^{y,H}$. 
	Suppose the end points of the spectral gap are both finite and denoted by $ x_j $, $ j = 1, 2 $. 
	Then, according to \eqref{7.3.2}, we have
	\begin{equation} 
	\begin{aligned}
	x_j=-\dfrac{1}{\ug_j}+y\int \dfrac{ t\mathrm{d} H(u,t)}{1 + u\ug_j + t\us_j},
	\\
	1=\dfrac{\us_j}{\ug_j}+y\int \dfrac{u\ug_j \mathrm{d} H(u,t)}{1 + u\ug_j + t\us_j}. \label{7.5.2}
	\end{aligned}
	\end{equation}
	Note that 
	\begin{align*}
	\frac{\partial x_j}{\partial \ug_j} =\frac{1}{\ug^2_j} - yA_{2} -yB_{2}\us'_j=0,
	\end{align*}
	Hence
	\begin{align}\label{7.5.3}
	\frac{\mathrm{d}x_j}{\mathrm{d}y}=\frac{\partial x_j}{\partial y}+\frac{\partial x_j}{\partial \ug}\frac{\mathrm{d} \ug}{\mathrm{d} y}=\int\dfrac{t\mathrm{d}H(u,t)}{1 + u\ug_j + t\us_j}.
	\end{align}
	Because $ (\ug_1, \us_1) $ and $ (\ug_2, \us_2) $ are in the same spectral gap, thus by the second
	conclusion of Lemma \ref{lemma2}, for any support point $ (u, t) $ of $ H $, 
	$$ (1 + u\ug_1 +t\us_1)(1 + u\ug_2 + t\us_2) > 0.  $$
	Therefore, we have
	\begin{align*}
	\frac{\mathrm{d}}{\mathrm{d}y}(x_1-x_2)=\int\dfrac{t[u(\ug_2-\ug_1)+t(\us_2-\us_1)]\mathrm{d}(u,t)}{(1 + u\ug_1 + t\us_1)(1 + u\ug_2 + t\us_2)}.
	\end{align*}
	
	By Lemma \ref{lemma2}, both $ x $ and $ \us $ are increasing functions of $ \ug $ outside the
	support of $F^{y,H}$ and hence $ (x_1-x_2) $ has the opposite signs as $ (\ug_2-\ug_1) $
	and $ (\us_2-\us_1) $, hence
	\begin{align*}
	\frac{\mathrm{d}}{\mathrm{d}y}\vert x_1-x_2\vert <0 ,
	\end{align*}
    thus, the size of spectral gap increases when $ y $ decreases.
    
    Finally, we consider the case where $ x_1 = \infty  $ and $ x_2  $ finite. 
    By \eqref{7.5.3}, we have
    \begin{align*}
    	\frac{\mathrm{d}x_2}{\mathrm{d}y}=\int\dfrac{t\mathrm{d}H(u,t)}{1 + u\ug_j + t\us_j}.
    \end{align*}
	Since $ (\ug, \us) $ is in the same spectral gap as $ (0, 0) $, we know that $ 1 + u\ug + t\us >
	0 $ and hence $ \mathrm{d}x_2/\mathrm{d}y < 0  $ which implies that $ x_2 $ decreases as $ y $ decreases.
	Therefore, in all cases, when $ y $ deceases, the size of the spectral gap increases.
	
	By Lemma \ref{lemma2}, there are $ \ug_a $ and $ \ug_b $ in $ \mathbb{R}\setminus \{0\}  $ such that $ a = x_{y,H}(\ug_a) $ and $ b = x_{y,H}(\ug_b) $. 
	Then, for any $ \ell $, we may select $ [a_{\ell}, b_{\ell}] $ such that $ b_{\ell}- a_{\ell} $ increases as $\ell$ increases. 
	By Lemma \ref{lemma5}, we have
	\begin{align} \label{lemma5.2.9}
		&\max_{k\leq p} \vert \lambda^{1/2}_k (\bbB_{n,\ell+1}) - \lambda^{1/2}_k (\bbB_{n,\ell}) \vert  \leq \| n_{\ell}^{-1/2} \bbT^{1/2} \bbX_{n_\ell} - n_{\ell+1}^{-1/2} \bbT^{1/2}\bbX_{n_{\ell+1}} \| \nonumber \\
\le&(n_{\ell}^{-1/2}-n_{\ell+1}^{-1/2})\|\bbT_n^{1/2}\bbX_{n_\ell}\|+n_{\ell+1}\|\bbX_\ell-\bbX_{n_{\ell+1}}\|\leq K\tau
	\end{align}
	which can be made arbitrarily small and hence we may assume that 
	$ a_{\ell+1} \in ( a_{\ell} - \frac{1}{2}(b_{\ell} - a_{\ell}), \frac{1}{2} (b_{\ell} + a_{\ell})  )$ and 
	$ b_{\ell+1} \in ( \frac{1}{2} (b_{\ell} + a_{\ell}), b_{\ell} + \frac{1}{2} (b_{\ell} - a_{\ell}),   ) $. 
	
\subsection{Proof of Theorem \ref{th731}}

    In Lemma \ref{lemma1}, it has proved that with probability one, when $ n $ is large, there are no eigenvalues of $ \bbB_n $ falling into the interval $ [a, b] $. 
    By the discussion in last section, there are intervals $ [a_\ell, b_\ell] $ satisfies the conditions of Lemma \ref{lemma1} such that when $ n $ is large, there are no eigenvalues of $ \bbB_{n,\ell} $ falling into the intervals $ [a_\ell, b_\ell] $.
    
    We claim that for each $ j \leq p $, $ \lambda_j(\bbB_{n,\ell})$  and  $ \lambda_j(\bbB_{n,\ell+1})$ are on the same sides of $ [a_\ell, b_\ell] $ and $ [a_{\ell+1}, b_{\ell+1}] $, respectively. 
    If it is not the case, then we should have for some $ j $ and $ \ell $,
    $$\vert\lambda_j(\bbB_{n,\ell+1})-\lambda_j(\bbB_{n,\ell})\vert \geq b_{\ell+1}\wedge b_\ell-a_{\ell+1} \vee a_\ell \geq \frac{1}{2}(b-a)  $$
    which contradicts \eqref{lemma5.2.9}.

	Therefore, the numbers of eigenvalues of $ \bbB_n $ on each side of the interval $[a, b]$ are exactly the same as those of eigenvalues of $\bbB_{n,\ell}$ on the corresponding sides of $[a_\ell, b_\ell]$ for any finitely $ \ell $. 
	When $ \ell\to\infty $, the matrix $ \bbB_{n,\ell} $ tends to the non-random $ \bbB_{n,\infty} = \bbR_n\bbR_n^* + \bbT_n $ whose paired eigenvalues are $ \{u_{n,i} + t_{n,i} \} $.

	The spectral gaps of the non-random $\bbB_{n,\infty}$ are the spacings of eigenvalues $u_{n,i}+t_{n,i}$. 
	And the limiting spectral gap must be a subinterval of some limiting spacings. 
	Suppose that the limiting spectral gap $ (a_\infty, b_\infty) $ is a subinterval of
	$(u_{i+1,n} +t_{i+1,n}, u_{i,n} +t_{i,n}) $ which satisfies $ d-c < u_{i,n} +t_{i,n}-(u_{i+1,n}+t_{i+1,n})$.
	By Lemma\ref{lemma5}, we have
	\begin{align}\label{7.5.4}
		&\max_{k\leq p} \vert \lambda_k^{1/2}(\bbB_{n,\ell})- \lambda_k^{1/2}(\bbB_{n,\infty})\vert \leq  \|n_\ell^{-1/2} \bbT_n^{1/2}\bbX_{n_\ell} -  \bbT_n^{1/2} \| \\ \nonumber
		\leq & K \max_{k\leq p} \| s_k(n^{-1}_{n_\ell} \bbX_{n_\ell} \bbB^*_{n,\ell}-\bbI_p)\| 
		\leq 3K \sqrt {y_\ell},
	\end{align}
	with probability $o(n_\ell^{-\mu})$ for any fixed $\mu>0$. Therefore, when $\ell$ is large so that $3K\sqrt{y_\ell}<\max(d-b,c-a)$, the numbers of eigenvalues of $\bbB_{n,\ell} $ larger than $ b_\ell $ (smaller than $ a_\ell$) are the same of eigenvalues of $\bbB_{n,\infty}$ larger than
	$u_{i,n}+t_{i,n}$ (smaller than $u_{i+1,n} +t_{i+1,n}$). 
	
	Now, we complete the proof of Theorem \ref{th731}. Note that for the non-random matrix $\bbB_{n,\infty}$, its eigenvalues are $ \{u_{i,n} + t_{i,n} \}$ and $ y_\infty = 0 $ which implies that $ \us(x) = \ug(x) = -1/x $.
	Therefore, $  h_j(x) = u_{j,n}\ug + t_{j,n}\us = -(u_{j,n} + t_{j,n})/x.  $ 
	Then, that $ h_j(x) $ is less than $ -1 $ is equivalent to $ u_{i,n} + t_{i,n} $ is larger than $ x $. Hence, by the continuity and monotonicity, the number eigenvalues $u_{i,n}+t_{i,n}$ are larger than 
	$ b_\infty $, is the same that $ h_j(x) $ less than $ -1 $ and that the eigenvalues are smaller than 
	$ a_\infty $ is the same that $h_j(x) > -1 $.
	
	In viewing \eqref{lemma5.2.9}, we may select $ \tau $ small enough so that 
	$$  \max_{k\leq p} \vert\lambda_k(\bbB_{n,{\ell+1}})-\lambda_k(\bbB_{n,\ell})\vert < (b - a)/2.  $$ 
	Also, when $ x \in [a_\ell, b_\ell] \cap [a_{\ell+1}, b_{\ell+1}] $, 
	$$  \vert \ug^{n,\ell+1}(x) - \ug^{n,\ell}(x)\vert < \max_{y-x>(b-a)/2}(\ug^{n,\ell}(y) - \ug^{n,\ell}(x))  $$ and
	$$   \vert \us^{n,\ell+1}(x) - \us^{n,\ell}(x)\vert < \max_{y-x>(b-a)/2}(\us^{n,\ell}(y) - \us^{n,\ell}(x)),  $$ 
	where $ \ug^{n,\ell} $ and $ \us^{n,\ell} $ are the Stieltjes transforms $ \ug $ and $ \us $ for the LSD of $ \bbB_{n,\ell} $ respectively. 
	Hence, if $ h^\ell_j(x) > -1 $ for all $ x\in [a_\ell, b_\ell] $, then $ h_j(a_\ell) > -1 $ and thus
	\begin{align*}
		&h^{\ell+1}_j(\frac{1}{2}(b_\ell+a_\ell)) = u_{jn}\ug^\ell(\frac{1}{2}(b_\ell+a_\ell)) + t_{jn} \us_\ell(\frac{1}{2}(b_\ell+a_\ell))\\
	\leq& u_{jn}\ug^\ell(a_\ell) + t_{jn} \us_\ell(a_\ell) = h^\ell_j(a_\ell) > -1.  
	\end{align*}


	Then, by the continuity and monotonicity of $h^{\ell+1}_j$, we conclude that
	$ h^{\ell+1}_j(x) > -1 $ for all $ x \in [a_{\ell+1}, b_{\ell+1}] $.

	Similarly, we can prove that $ h^{\ell+1}_j(x) < -1 $ for all $ x \in [a_{\ell+1}, b_{\ell+1}] $ if $ h^\ell_j(x) < -1 $ for $ x\in[a_\ell, b_\ell] $.
	
	Finally, by \eqref{7.5.4}, we can prove that for the last $ \ell $, $ h^\ell_j(x) $ is larger (smaller)
	than $ -1 $ are the same as $ \bbB_{n,\infty} $. Thus, Theorem \ref{th731} is proved.


	\begin{acks}[Acknowledgments]
		The authors would like to thank the anonymous referee, the Associate Editor and the Editor for their invaluable and constructive comments.
		J. Hu was supported by NSFC (No. 12171078, 11971097) and National Key R \& D Program of China (No. 2020YFA0714102).
		Z. D. Bai was partially supported by NSFC Grant 12171198 and 12271536 and Team Project of Jilin Provincial Department of Science and Technology (No.20210101147JC). 
		
	\end{acks}

	\bibliographystyle{imsart-number} 

\end{document}